\documentclass{article}

\usepackage{arxiv}

\usepackage{amsmath,amssymb}

\usepackage[utf8]{inputenc} 
\usepackage[T1]{fontenc}    
\usepackage{hyperref}       
\usepackage{url}            
\usepackage{booktabs}       
\usepackage{amsfonts}       
\usepackage{nicefrac}       
\usepackage{microtype}      
\usepackage{cleveref}       
\usepackage{lipsum}         
\usepackage{graphicx}
\usepackage{natbib}
\usepackage{doi}

\usepackage{listings}
\usepackage{color}
\usepackage{caption}
\usepackage[labelformat=simple]{subcaption}

\usepackage{multirow}

\usepackage{cprotect}

\newcommand{\vect}[1]{\boldsymbol{#1}}
\newcommand{\mat}[1]{\boldsymbol{#1}}

\definecolor{codegreen}{rgb}{0,0.6,0}
\definecolor{codegray}{rgb}{0.5,0.5,0.5}
\definecolor{codepurple}{rgb}{0.58,0,0.82}
\definecolor{backcolour}{rgb}{0.95,0.95,0.92}

\lstdefinestyle{mystyle}{
    commentstyle=\color{codegreen},
    keywordstyle=\color{magenta},
    numberstyle=\tiny\color{codegray},
    stringstyle=\color{codepurple},
    basicstyle=\color{black}\ttfamily\normalsize,
    breakatwhitespace=false,
    breaklines=true,
    captionpos=b,
    keepspaces=true,
    showspaces=false,
    showstringspaces=false,
    showtabs=false,
    tabsize=1,
    frame=single,
    numbers=left,
    xleftmargin=18pt,
    xrightmargin=5pt
}
\usepackage{booktabs}
\usepackage{hyperref}

\usepackage{comment}

\usepackage{colortbl}
\usepackage{siunitx}
\newcommand{\method}{F3R}
\newcommand{\methodD}{fp64-\method{}}
\newcommand{\methodF}{fp32-\method{}}
\newcommand{\methodH}{fp16-\method{}}
\newcommand{\methodHBest}{fp16-\method{}-best}
\newcommand{\double}{fp64}
\newcommand{\single}{fp32}
\newcommand{\half}{fp16}
\newcommand{\dbar}[1]{\bar{\bar{#1}}}

\usepackage[linesnumbered,ruled]{algorithm2e}
\usepackage{enumitem}

\begin{document}

\title{A Nested Krylov Method \\Using Half-Precision Arithmetic}

\author{ \href{https://orcid.org/0009-0002-8328-0388}{\includegraphics[scale=0.06]{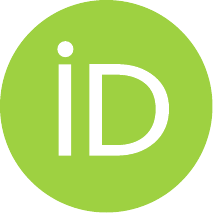}\hspace{1mm}Kengo~Suzuki} \\
    Academic Center for Computing and Media Studies, \\
    Kyoto University, \\
    Kyoto, Japan \\
    \texttt{suzuki@i.kyoto-u.ac.jp} \\
    \And
    Takeshi~Iwashita \\
    Academic Center for Computing and Media Studies, \\
    Kyoto University, \\
    Kyoto, Japan \\
    \texttt{iwashita@i.kyoto-u.ac.jp} \\
}

\renewcommand{\shorttitle}{A Nested Krylov Method Using Half-Precision Arithmetic}

\maketitle

\begin{abstract}
Low-precision computing is essential for efficiently utilizing memory bandwidth and computing cores.
While many mixed-precision algorithms have been developed for iterative sparse linear solvers, effectively leveraging half-precision (\half{}) arithmetic remains challenging.
This study introduces a novel nested Krylov approach that integrates the FGMRES and Richardson methods in a deeply nested structure, progressively reducing precision from double-precision to \half{} toward the innermost solver.
To avoid meaningless computations beyond precision limits, the low-precision inner solvers perform only a few iterations per invocation, while the nested structure ensures their frequent execution.
Numerical experiments show that incorporating \half{} into the approach directly enhances solver performance without compromising convergence, achieving speedups of up to $1.65\times$ and $2.42\times$ over double-precision and double-single mixed-precision implementations, respectively.
Furthermore, the proposed method outperforms conventional mixed-precision Krylov solvers, CG, BiCGStab, and restarted FGMRES, by factors of up to 2.47, 2.74, and 69.10, respectively.
\end{abstract}


\section{Introduction} \label{sec:intro}
We explore solving sparse linear systems, $\mat{A}\vect{x} = \vect{b}$, where $\mat{A} \in \mathbb{R}^{n \times n}$ and $\vect{b}, \vect{x} \in \mathbb{R}^{n}$, using iterative methods with low-precision arithmetic and data types.
Because most computing kernels in sparse iterative methods are memory-bounded, reducing data movement directly improves performance.
Using low-precision arithmetic is also motivated by a recent trend in hardware development where computing units are simplified to increase parallelism.
Hence, significant efforts have been made to realize practical sparse linear solvers in low-precision, with many proposals and evaluations~\citep{abdelfattah2021survey,higham2022mixed}.
Among these, mixed-precision algorithms, particularly Generalized Minimal Residual (GMRES) algorithms, have received considerable attention~\citep{turner1992efficient,anzt2011error,lindquist2021accelerating,zhao2022numerical,carson2022mixed}, leading to a proposal of a mixed-precision benchmark, HPGMP~\citep{yamazaki2022high}, as a successor to the HPCG benchmark.
While their practicality has steadily improved, there is room for improvement in using half-precision (\half{}).
In particular, despite its prevalence in direct methods, \half{} remains relatively unexplored in sparse iterative methods.
Furthermore, existing research mainly focuses on the precision of preconditioners~\citep{carson2023mixed,scott2024avoiding}, with very few studies investigating the use of \half{} in the main part of sparse iterative methods.
This rareness is primarily due to the significant degradation in the convergence speed by \half{} even in algorithms where single-precision (\single{}) works well~\citep{aliaga2023compressed}.

To leverage \half{} arithmetic effectively, we introduce a novel nested Krylov-based approach that integrates flexible GMRES (FGMRES)~\citep{saad2003iterative} and Richardson~\citep{richardson1911ix} methods.
As shown in previous work (e.g., Reference~\citep{yamazaki2022high}), low-precision GMRES can partially reduce residuals, even when using a unified precision (e.g., \single{} or \half{}).
Our approach exploits this property by progressively reducing both precision and the number of iterations toward the innermost solver.
Low-precision inner solvers perform only a few iterations per invocation to prevent exceeding precision limits, while the nested structure increases the proportion of low-precision computations.
We also provide an adaptive technique for updating the weight in the \half{} Richardson part, which is located in the innermost and affects solver performance significantly.
The adaptive technique enhances the stability of our method.
We perform numerical experiments on both CPU and GPU nodes and show that our method outperforms standard Krylov methods, including restarted FGMRES, Conjugate Gradient (CG) and Bi-CG Stabilized (BiCGStab)~\citep{saad2003iterative}.

The rest of this paper is structured as follows.
In Section~\ref{sec:related}, we review previous mixed-precision algorithms and the nested Krylov framework.
Here, we also clarify the difference between this study and previous studies.
In Section~\ref{sec:notation}, we introduce a notation and terminology for nested Krylov.
In Section~\ref{sec:propose}, we explain the proposed approach.
First, we evaluate it in comparison with the restarted FGMRES, CG, and BiCGStab methods in Section~\ref{sec:numerical}.
Next, we give more detailed analyses of its parameter setting in Section~\ref{sec:ablation}.
Finally, we conclude the study in Section~\ref{sec:conclution}.

\section{Related Work and Contributions}
\label{sec:related}
Mixed-precision algorithms are widely used across various applications in numerical linear algebra~\citep{abdelfattah2021survey,higham2022mixed}, with mixed-precision solvers developed for both dense and sparse linear systems.
For solvers that incorporate direct methods like LU (Cholesky) factorization, numerous techniques exploit \single{} and \half{}~\citep{carson2017new, haidar2018harnessing, abdelfattah2020investigating, amestoy2024five}.

In sparse iterative solvers, which we focus on, low-precision arithmetic is typically integrated in two ways: within parts of the solver or by nesting a low-precision solver as a preconditioner within an outer high-precision solver.
The former approach often involves reducing the precision of preconditioners.
Many techniques have been proposed to utilize \single{} and \half{} in preconditioners, including incomplete Cholesky, approximate inverse, and multigrid preconditioners~\citep{zong2024fp16,scott2024avoiding,carson2023mixed}.
Another strategy is to reduce the precision of the sparse matrix-vector multiplication (SpMV) kernel, a core component of sparse iterative solvers.
Several studies have shown the viability of low-precision SpMV in Krylov solvers~\citep{simoncini2003theory, graillat2024adaptive}.
Additionally, some methods reduce the precision of specific variables in the solver, such as the Arnoldi basis in GMRES~\citep{aliaga2023compressed, grutzmacher2024frsz2}, where 32-bit types can be effective, though \half{} offers less benefit.

The latter approach frequently employs iterative refinement, where a low-precision solver serves as a preconditioner for a high-precision Richardson solver.
GMRES is commonly used because of its affinity with restarting included in iterative refinement.
The effectiveness of using GMRES implemented using \single{} or 32-bit integers as the inner solver has been explored both numerically and analytically~\citep{anzt2011error,loe2021experimental,zhao2022numerical,lindquist2021accelerating,suzuki2024integer}.
Using other iterative methods, such as CG and BiCGStab, as the inner solver have also been investigated~\citep{zhao2023numerical,higham2021exploiting}.
Furthermore, some studies uses more sophisticated solvers, such as double-precision (\double{}) FGMRES and CG, as the outer solver instead of iterative refinement~\citep{baboulin2009accelerating,buttari2008using}.

Another important technique related to this study is nested Krylov, a framework that recursively uses one Krylov solver as a preconditioner for another Krylov solver.
One of the most well-known instances of it would be the flexible inner-outer GMRES method~\citep{saad1993flexible}, where a preconditioned GMRES method is used as a preconditioner of the FGMRES method.
Besides, several other variations have been developed individually for specific purposes~\citep{golub1999inexact, axelsson2004class}.
The nested Krylov framework is discussed systematically in Reference~\citep{mcinnes2014hierarchical}, particularly within the context of reducing communication overhead in distributed systems.
Note that nested Krylov can, in principle, accept any iterative methods other than the Krylov method despite its name, as the Chebyshev iteration is used in Reference~\citep{mcinnes2014hierarchical}.

\subsection{Contributions}
While these previous studies on mixed-precision algorithms have mainly focused on two-level nesting and using \single{}, this study distinguishes itself by exploiting low-precision arithmetic in three- or higher-level nested structures and investigating the use of \half{}.
Consequently, the main contributions of this study are as follows:
\begin{itemize}[leftmargin=2em]
\item We propose a mixed-precision (\double{}, \single{}, and \half{}) approach by integrating FGMRES and Richardson in the nested Krylov framework. In the approach, we gradually reduce the precision towards the innermost solver to achieve effective use of \half{}.
\item We adopt an \half{} Richardson solver in the innermost part and also propose a technique for controlling its weight adaptively.
\item We show that our solver outperforms conventional mixed-precision solvers through numerical experiments on CPU and GPU nodes. Notably, our solver is superior to preconditioned CG even for several symmetric positive definite matrices.
\end{itemize}

\section{Notation} \label{sec:notation}
For clarity, we express nested Krylov solvers using tuples.
A pair $(\mat{S}^{(1)}, \mat{M})$ denotes a solver $\mat{S}^{(1)}$ preconditioned by $\mat{M}$.
Similarly, if $\mat{S}^{(1)}$ is preconditioned by another nested solver ($\mat{S}^{(2)}$, $\mat{M}$), it is written as ($\mat{S}^{(1)}$, ($\mat{S}^{(2)}$, $\mat{M}$)).
For simplicity, we also express it as ($\mat{S}^{(1)}$, $\mat{S}^{(2)}$, $\mat{M}$) following the standard notation for ordered tuples.
This naturally extends to deeper levels of nesting, resulting in the general form ($\mat{S}^{(1)}$, $\mat{S}^{(2)}$, $\ldots$, $\mat{S}^{(D)}$, $\mat{M}$), where $D$ is the depth level of the nested solver.
For readability, we introduce the following terminology.
The solvers $\mat{S}^{(2)}$ through $\mat{S}^{(D)}$ are referred to as inner solvers, while $\mat{S}^{(1)}$ and $\mat{S}^{(D)}$ are called the outermost and innermost solvers, respectively.
Moreover, $\mat{S}^{(d - 1)}$ (for $d = 2,\ldots,D$) is referred to as the parent solver of $\mat{S}^{(d)}$.
The preconditioner $\mat{M}$ is called the primary preconditioner because an inner solver $\mat{S}^{(d)}$ can also be regarded as a preconditioner of its parent solver.
In this study, $\mat{M}$ is assumed to be an algebraic preconditioner, such as an incomplete LU (ILU) preconditioner.
Lastly, FGMRES and Richardson of $m$ iterations are denoted as $\mat{F}_m$ and $\mat{R}_m$, respectively.

\section{Nested Krylov-Based Approach} \label{sec:propose}
Previous studies have shown that the performance of low-precision GMRES in the early iterations can be comparable to that of full-precision GMRES, even in the case of uniform-precision implementations~\citep{yamazaki2022high,zhao2022numerical}.
This suggests that low-precision arithmetic is effective when the number of iterations per invocation is limited.
Building on this insight we propose a novel mixed-precision approach that utilizes \half{} arithmetic effectively.
Our method employs a nested Krylov structure, integrating multiple FGMRES and Richardson solvers while gradually reducing precision from \double{} to \half{} and decreasing the number of iterations towards the innermost solver.
Since each inner solver runs only a few iterations per invocation, the negative effect of using low precision on convergence is expected to be small.
However, because these inner solvers are invoked at every iteration of their parent solvers, low-precision arithmetic is performed frequently.
As a result, our approach can achieve a high frequency of \half{} computations without deteriorating convergence.

\subsection{Our Strategy for Nesting}
\label{subsec:guide-to-nesting}
In our approach, we have various possible settings for the nesting depth and the number of iterations of each part.
To find effective configurations, we consider convergence speed and the amount of memory accesses per iteration, both of which are key factors in solver performance.
First, we make the following assumptions with respect to the convergence speed, whose validity is examined in Section~\ref{sec:ablation}:
\begin{enumerate}[label=(\roman*),leftmargin=2em]
\item \label{enum:fgmres} when $m$ is relatively large, for $\bar{m}$ and $\dbar{m}$ such that $m = \bar{m}\dbar{m}$, we can approximate ($\mat{F}_m$, $\mat{M}$) with a nested solver ($\mat{F}_{\bar{m}}$, $\mat{F}_{\dbar{m}}$, $\mat{M}$) without reducing the convergence speed.
\item \label{enum:richardson} when $m$ is small (e.g., $m < 5$), we can replace $(\mat{F}_{m}, \mat{M})$ with $(\mat{R}_{m}, \mat{M})$ without compromising the convergence speed.
\end{enumerate}
We expect these assumptions to hold when the primary preconditioner $\mat{M}$ is sufficiently effective.

Under these assumptions, we focus on the amount of memory accesses, the other key factor to achieving high performance.
We therefore introduce a rough model of the amount of memory accesses per $n$ of $(\mat{F}_{m}, \mat{M})$, $O_{(\mat{F}_{m}, \mat{M})}$, and that of $(\mat{R}_{m}, \mat{M})$, $O_{(\mat{R}_{m}, \mat{M})}$, as
\begin{equation}
\begin{split}
&O_{(\mat{F}_{m}, \mat{M})} := c_Am + c_Mm + \frac{5}{2}m^{2}, \\ 
&O_{(\mat{R}_{m}, \mat{M})} := c_A(m-1) + c_Mm + 4(m-1),
\end{split}
\end{equation}
where $c_A$ and $c_M$ are constants associated with the number of nonzeros per row of the coefficient matrix $\mat{A}$ and preconditioner $\mat{M}$, respectively.
Terms $c_A(m-1)$ and $4(m-1)$ in $O_{(\mat{R}_{m}, \mat{M})}$ reflect that we assume the initial guess to be a zero vector.
Note that these models are intended to provide a rough guideline, and we here do not consider exact models, which must include complicated factors such as data types, cache utilization, and implementation details.

Then, we consider the memory access model of two-level nested FGMRES ($\mat{F}_{\bar{m}}$, $\mat{F}_{\dbar{m}}$, $\mat{M}$), $O_{(\mat{F}_{\bar{m}}, \mat{F}_{\dbar{m}}, \mat{M})}$, while assuming $m = \bar{m}\dbar{m}$ to fix the number of calls to the primary preconditioner $\mat{M}$:
\begin{equation}{}\label{eq::fgmres-guide}
\begin{split}
O_{(\mat{F}_{\bar{m}}, \mat{F}_{\dbar{m}}, \mat{M})} &:= c_A\bar{m} + O_{(\mat{F}_{\dbar{m}}, \mat{M})}\bar{m} + \frac{5}{2}\bar{m}^{2} \\
& = O_{(\mat{F}_m, \mat{M})} + c_A\bar{m} + \frac{5}{2}\dbar{m}^{2}\bar{m} + \frac{5}{2}\bar{m}^{2} - \frac{5}{2}m^{2}.
\end{split}
\end{equation}
This equation provides an estimation of the overhead and benefits of dividing reference FGMRES into two-level nested FGMRES.
When $m$ is relatively large, nesting offers a significant advantage due to the reduced cost of the Arnoldi process.
In other words, in such cases, splitting FGMRES has the benefit of both reducing memory accesses and improving the usability of low precision by reducing the number of iterations per call.
For example, assuming $c_A = 45$ (30 nonzeros per row, with \double{} for values and 32-bit integers for indices) and $m = 64$, the inequity $O_{(\mat{F}_{\bar{m}}, \mat{F}_{\dbar{m}}, \mat{M})}$ $<$ $O_{(\mat{F}_m, \mat{M})}$ holds for most possible values of $\bar{m}$, and $\bar{m} = 10$ results in the least amount, though 10 is not a divisor of 64.
On the other hand, for small values of $m$, Equation~\eqref{eq::fgmres-guide} also indicates that nesting rather increases memory accesses.
This would be a problem when $m$ is still too large to implement the inner solver in lower-precision such as \half{}.
However, if we replace the inner FGMRES with Richardson, the amount becomes
\begin{equation}{}\label{eq::richardson-guide}
\begin{split}
O_{(\mat{F}_{\bar{m}}, \mat{R}_{\dbar{m}}, \mat{M})} &:= c_A\bar{m} + O_{(\mat{R}_{\dbar{m}}, \mat{M})}\bar{m} + \frac{5}{2}\bar{m}^{2} \\
& = O_{(\mat{F}_m, \mat{M})} + 4(\dbar{m}-1)\bar{m} + \frac{5}{2}\bar{m}^{2} - \frac{5}{2}m^{2}.
\end{split}
\end{equation}
This leads to fewer memory accesses for all possible values of $\bar{m}$ when $m \geq 3$.
Thus, we will achieve higher performance by replacing the inner FGMRES with Richardson as long as Assumption~\ref{enum:richardson} remains valid.

In summary, we split a baseline FGMRES solver into a nested solver, where inner solvers perform a few iterations for the use of low-precision arithmetic.
To achieve this, we take the following strategy.
Building on Assumption~\ref{enum:fgmres}, we recursively divide an inner FGMRES solver into a nested FGMRES solver while considering reducing memory accesses based on Equation~\eqref{eq::fgmres-guide}.
Even when Equation~\eqref{eq::fgmres-guide} indicates that nesting increases memory accesses, replacing FGMRES with Richardson can reduce memory accesses as shown in Equation~\eqref{eq::richardson-guide}.
Therefore, while Assumption~\ref{enum:richardson} holds, we may consider using Richardson instead of inner FGMRES.

\subsection{Proposed Solver: \method{}} \label{subsec:fxr}
Following the strategy outlined above, we develop an instance of our nested Krylov-based approach, named \method{}.
It consists of three FGMRES solvers and an innermost Richardson solver, defined as ($\mat{F}_{m_1}$, $\mat{F}_{m_2}$, $\mat{F}_{m_3}$, $\mat{R}_{m_4}$, $\mat{M}$) in the tuple notation, with default parameters ($m_1$, $m_2$, $m_3$, $m_4$) = (100, 8, 4, 2).
The development process of \method{} is as follows.

First, we select ($\mat{F}_{64}$, $\mat{M}$) as a reference.
Subsequently, we approximate it with ($\mat{F}_{m_2}$, $\mat{F}_{m_3'}$, $\mat{M}$).
Given the considerations of the amount of memory accesses in Equation~\eqref{eq::fgmres-guide}, we choose $m_2 = 8$, a value considered near-optimal among the divisors of 64; then, $m_3' = 8 (= 64/m_2)$.
Presuming that $m_3' = 8$ is still too large for implementing FGMRES in \half{}, we attempt to further divide the inner solver ($\mat{F}_{m_3'}$, $\mat{M}$) into ($\mat{F}_{m_3}$, $\mat{F}_{m_4}$, $\mat{M}$).
However, this will increase memory accesses, according to Equation~\eqref{eq::fgmres-guide}.
Thus, we replace $\mat{F}_{m_4}$ with $\mat{R}_{m_4}$;
that is, we divide ($\mat{F}_{m_3'}$, $\mat{M}$) into ($\mat{F}_{m_3}$, $\mat{R}_{m_4}$, $\mat{M}$).
While there are several possible choices for $m_4$, we opt $m_4 = 2$ considering Assumption~\ref{enum:richardson}.
Finally, we wrap ($\mat{F}_{m_2}$, $\mat{F}_{m_3}$, $\mat{R}_{m_4}$, $\mat{M}$) with $\mat{F}_{m_1}$, rather than restarting (or iterative refinement), for better performance, and we set $m_1 = 100$ empirically.

Convergence is checked only in the outermost $\mat{F}_{m_1}$.
In practice, if convergence is not reached within 100 outermost iterations, we execute the entire solver ($\mat{F}_{m_1}$, $\mat{F}_{m_2}$, $\mat{F}_{m_3}$, $\mat{R}_{m_4}$, $\mat{M}$) again in the manner of the restarting technique.

We implement \method{} in mixed precision as specified in Table~\ref{tab:precision}.
Starting from \double{} for the outermost FGMRES $\mat{F}_{m_1}$, we progressively reduce precision towards the innermost \half{} Richardson ($\mat{R}_{m_4}$, $\mat{M}$).
The Richardson method has a straightforward recurrence and avoids reduction kernels, such as the norm kernel, which are susceptible to numerical errors.
Moreover, we set $m_4$ to 2, a small value.
Thus, using \half{} would not significantly affect computational accuracy.
The two inner FGMRES solvers, $\mat{F}_{m_2}$ and $\mat{F}_{m_3}$, are also implemented in low precision.
These solvers include reduction kernels, so we store their variables in \single{}.
Specifically, as we employ classical Gram-Schmidt and Givens rotation for the Arnoldi process and the QR factorization, all associated computations are performed only with vectors and scalars stored in \single{}.
While $\mat{F}_{m_2}$ employs an \single{} SpMV kernel, $\mat{F}_{m_3}$ stores $\mat{A}$ in \half{}. 
This choice is based on the fact that $\mat{F}_{m_3}$ serves as a preconditioner for $\mat{F}_{m_2}$, where higher precision would not be necessitated.
Finally, $\mat{F}_{m_1}$ is implemented in \double{} to provide an accurate approximate solution.
Throughout the implementation, higher-precision instructions are used when the inputs differ in precision.
For example, $\mat{F}_{m_3}$ performs SpMV in \single{} because $\mat{A}$ is stored in \half{} while the input Arnoldi basis is in \single{}.

One may assume that the setting of $m_4 = 2$ for the innermost \half{} Richardson solver is too small.
However, we emphasize that the innermost solver performs $m_2m_3m_4$ iterations per outermost iteration in total, achieving a high frequency of \half{} computations.

\begin{table}[tb]
\caption{Details of the proposed method \method{}.}
\label{tab:precision}
\small
\centering
\begin{tabular}{llrrr}
\toprule
         & & \multicolumn{3}{c}{Precision} \\\cmidrule{3-5}
Solver   & Default & $\mat{A}$       & Vectors & $\mat{M}$ \\
\midrule
$\mat{F}_{m_1}$ & $m_1 = 100$ & \double{} & \double{}  & -   \\
$\mat{F}_{m_2}$ & $m_2 = 8$ & \single{} & \single{}  & -   \\
$\mat{F}_{m_3}$ & $m_3 = 4$ & \half{}   & \single{}  & -   \\
$\mat{R}_{m_4}$ & $m_4 = 2$ & \half{}   & \half{}    & \half{} \\
\bottomrule
\end{tabular}
\end{table}

\begin{algorithm}[tb]
\caption{Richardson with adaptive weight-updating}
\label{alg:richardson}
\DontPrintSemicolon
\KwData{Linear system $\mat{A}\vect{z} = \vect{v}$, preconditioner $\mat{M}$, weights $\omega_1,\ldots,\omega_{m_4}$, and call count of Algorithm~\ref{alg:richardson} $c_{ntr}$.}
\KwResult{Approximate solution $\vect{z}_{m_4}$, weights $\omega_1,\ldots,\omega_{m_4}$, and call count of this algorithm $c_{ntr}$.}
\tcc{$\omega_1,\ldots,\omega_{m_4}$ and $c_{ntr}$ are global variables and retain their values across function calls.}
\If{$c_{ntr} = 1$ (the first call)}{
  Initialize $\omega_1 = \cdots = \omega_{m_4} = 1$. \;
}
Assume $\vect{z}_{0} = \vect{0}$. \;
\For{$k := 1,\ldots, m_4$}{ 
  $\vect{r}_{k-1} := \vect{v} - \mat{A}\vect{z}_{k-1}$ \tcp*{$\vect{r}_{0} = \vect{v}$ without calculation.}
  \eIf {$c_{ntr} \equiv 0 \pmod{c}$} {
    $\omega'_k = (\vect{r}_{k-1}, \mat{A}\mat{M}\vect{r}_{k-1}) / (\mat{A}\mat{M}\vect{r}_{k-1}, \mat{A}\mat{M}\vect{r}_{k-1})$ \;
    $\vect{z}_{k} = \vect{z}_{k-1} + \omega'_k \mat{M}\vect{r}_{k-1}$ \;
    $\omega_k := (l\omega_k + \omega'_k)/(l + 1)$ for $l$ := $c_{ntr}$ / $c$\;
  } {
    $\vect{z}_{k} = \vect{z}_{k-1} + \omega_k \mat{M}\vect{r}_{k-1}$ \;
  }
}
Increment $c_{ntr} := c_{ntr} + 1$\;
\end{algorithm}

\subsection{Adaptive Weight Updating}
Finally, we consider the weight in the innermost Richardson part, ($\mat{R}_{m_4}$, $\mat{M}$), of \method{}, which is crucial for holding Assumption~\ref{enum:richardson}.
This Richardson part receives a vector $\vect{v}$ and searches an approximate solution of a linear system $\mat{A}\vect{z} = \vect{v}$ as a preconditioning step of its parent solver.
Specifically, at the $k$-th ($k = 1,\ldots,m_4$) iteration, it computes an approximate solution $\vect{z}_{k}$ as 
\begin{equation}
\begin{split}
\vect{z}_{k} &= \vect{z}_{k-1} + \omega \mat{M}(\vect{v} - \mat{A}\vect{z}_{k-1}),
\end{split}
\end{equation}
where $\omega$ is a scalar value, weight. 
Since Richardson is a stationary iterative method, its convergence speed depends on the spectral radius of $\mat{I} - \omega\mat{M}\mat{A}$.
Therefore, to satisfy Assumption~\ref{enum:richardson}, it is important to consider not only the choice of the primary preconditioner $\mat{M}$ but also the quality of the weight $\omega$.
Although the optimal weight depends on the eigenvalues of $\mat{M}\mat{A}$ and is impractical to compute, a locally optimal value can be easily determined at each iteration.
The $k$-th residual $\vect{r}_{k}$ in Richardson is given as
$$
\vect{r}_{k} = \vect{r}_{k-1} - \omega \mat{A}\mat{M}\vect{r}_{k-1}.
$$
Thus, the weight $\omega'_k$ that minimizes the 2-norm of $\vect{r}_{k}$ is found as 
$$
\omega'_{k} = \frac{(\vect{r}_{k-1}, \mat{A}\mat{M}\vect{r}_{k-1})}{(\mat{A}\mat{M}\vect{r}_{k-1}, \mat{A}\mat{M}\vect{r}_{k-1})}.
$$
If $\omega'_{k}$ is computed at every iteration, it is equivalent to GMRES(1), and Assumption~\ref{enum:richardson} may have some validity.
However, this approach undermines the advantage of Richardson, as it incurs additional SpMV and reduction computations.
Conversely, using a single $\omega'_{k}$ persistently throughout the whole iteration or for the entire solution process may also be ineffective;
there is no guarantee that $\omega'_{k}$ remains optimal in subsequent iterations or for different linear systems at other preconditioning steps of the parent solver.

Therefore, we propose a technique to estimate an appropriate value $\omega_k$ for each iteration by updating it over successive invocations of the Richardson solver.
While the optimal value of $\omega$ is unique, we use $k$ different weights $\omega_k$ for flexibility.
Specifically, we initialize $\omega_k$ as 1 and update it by computing $\omega'_{k}$ every $c$ calls to the Richardson solver, using the cumulative average of all previously computed weights:
\begin{equation}
\omega_k := \frac{l \omega_k + \omega'_k}{l + 1},
\end{equation}
where $l$ is the update count, defined as the number of invocations of the Richardson solver divided by $c$.
Importantly, this adaptive technique updates $\omega_k$ across all invocations of the Richardson part because the optimal weight does not depend on the right-hand side of a given linear system.
Consequently, the technique is expected to prevent $\omega_k$ from becoming overly localized, while keeping the additional computational cost low.
Algorithm~\ref{alg:richardson} outlines the complete process of the Richardson part with this adaptive weight-updating technique.
It is noteworthy that when updating $\omega_k$ on lines~7--11, $\omega'_{k}$ is used as the weight, rather than $\omega_{k}$, because it minimizes the residual at that step.
The effectiveness of this technique is evaluated through numerical experiments in Section~\ref{sec:ablation}.
It should be mentioned that the process for computing $\omega'_{k}$ is implemented in \single{}, although other parts of the innermost Richardson is implemented in \half{} as explained in Section~\ref{subsec:fxr}.

Our nested Krylov approach, in principle, accepts Richardson at multiple levels.
Note that in such cases, the weight $\omega_k$ should be handled separately at each level.
This is because the optimal value of a Richardson part depends on its preconditioner (inner solver), not the primary preconditioner, which differs from level to level.

\section{Experimental Results}
\label{sec:numerical}
In this section, we evaluate the effect of reducing precision on the performance of \method{} and compare \method{} with three conventional preconditioned Krylov methods, restarted FGMRES of a restart cycle of 64 (FGMRES(64)), CG, and BiCGStab, both on CPU and GPU nodes.
In addition to the setting outlined in Table~\ref{tab:precision}, we implemented \method{} in only \double{} or \double{}-\single{} mixed precision; the latter use \single{} for all the inner solvers.
To distinguish them, we refer to the implementation that uses three precisions, \double{}, \single{}, and \half{}, shown in Table~\ref{tab:precision}, as \methodH{}, while we label the additional two as \methodD{} and \methodF{}, respectively.
Leaving the detailed examination of the parameters, $m_1$, $m_2$, $m_3$, $m_4$, and $c$, to Section~\ref{sec:ablation}, this section focuses primarily on results with the default setting ($m_1 = 100$, $m_2 = 8$, $m_3 = 4$, $m_4 = 2$, and $c = 64$).
We also present results for \methodH{} configured for optimal performance, which is referred to as \methodHBest{}.
For each of the other Krylov methods, FGMRES(64), CG, and BiCGStab, we developed three \double{}-based solvers, varying the precision of the preconditioner storage between \double{}, \single{}, and \half{}.
When reducing the precision of the preconditioner, we first construct it in \double{} and then cast its values to \single{} or \half{}.
Similar to \method{}, these are prefixed by \double{}, \single{}, and \half{}: e.g., fp16-CG for \double{} CG with an \half{} preconditioner.
All the solvers used 32-bit integers for column indices and index pointer arrays in storing matrices.

\begin{table}[tb]
\small
\caption{Test matrices. $n_{nz}$ denotes the number of nonzeros. 
The suffix \texttt{x\_y\_z} for matrices from HPCG and HPGMP denotes the problem size in log2 for each of the x, y, and z axes.}
\centering
\begin{tabular}{lS[table-format=8.0]S[table-format=9.0]rrrr}
\toprule
Matrix & \multicolumn{1}{r}{$n$} & \multicolumn{1}{r}{$n_{nz}$} & $n_{nz}$/$n$ & $\alpha_\mathrm{ILU}$ & $\alpha_\mathrm{AINV}$ \\
\midrule
\texttt{Bump\_2911} & 2911419 & 127729899 & 43.87 & 1.1 & 1.2 \\
\texttt{Emilia\_923} & 923136  &  40373538 & 43.74 & 1.0 & 1.2 \\
\texttt{G3\_circuit} & 1585478 &   7660826 & 4.83  & 1.0 & 1.0 \\
\texttt{Queen\_4147} & 4147110 & 316548962 & 76.33 & 1.1 & 1.3 \\
\texttt{Serena} & 1391349 &  64131971 & 46.09 & 1.1 & 1.2 \\
\texttt{apache2} & 715176  &   4817870 & 6.74  & 1.0 & 1.0 \\
\texttt{audikw\_1} & 943695  &  77651847 & 82.28 & 1.1 & 1.6 \\
\texttt{ecology2} & 999999  &   4995991 & 5.00  & 1.0 & 1.0 \\
\texttt{hpcg\_7\_7\_7} & 2097152 &  55742968 & 26.58 & 1.0 & 1.0 \\
\texttt{hpcg\_8\_7\_7} & 4194304 & 111777784 & 26.65 & 1.0 & 1.0 \\
\texttt{hpcg\_8\_8\_7} & 8388608 & 224140792 & 26.72 & 1.0 & 1.0 \\
\texttt{hpcg\_8\_8\_8} &16777216 & 449455096 & 26.79 & 1.0 & 1.0 \\
\texttt{ldoor} & 952203  &  42493817 & 44.63 & 1.1 & 1.3 \\
\texttt{thermal2} & 1228045 &   8580313 & 6.99  & 1.0 & 1.0 \\
\texttt{tmt\_sym} & 726713  &   5080961 & 6.99  & 1.0 & 1.0 \\
\midrule
\texttt{Freescale1     }& 3428755  &  17052626 & 4.97  & 1.1 & 1.1 \\
\texttt{Transport      }& 1602111  &  23487281 & 14.66 & 1.0 & 1.0 \\
\texttt{atmosmodd      }& 1270432  &   8814880 & 6.94  & 1.0 & 1.0 \\
\texttt{atmosmodj      }& 1270432  &   8814880 & 6.94  & 1.0 & 1.0 \\
\texttt{atmosmodl      }& 1489752  &  10319760 & 6.93  & 1.0 & 1.0 \\
\texttt{hpgmp\_7\_7\_7 }&  2097152 &  55742968 & 26.58 & 1.0 & 1.0 \\
\texttt{hpgmp\_8\_7\_7 }&  4194304 & 111777784 & 26.65 & 1.0 & 1.0 \\
\texttt{hpgmp\_8\_8\_7 }&  8388608 & 224140792 & 26.72 & 1.0 & 1.0 \\
\texttt{hpgmp\_8\_8\_8 }& 16777216 & 449455096 & 26.79 & 1.0 & 1.0 \\
\texttt{rajat31        }& 4690002  &  20316253 & 4.33  & 1.0 & 1.0 \\
\texttt{ss             }& 1652680  &  34753577 & 21.03 & 1.1 & 1.2 \\
\texttt{stokes         }& 11449533 & 349321980 & 30.51 & 1.0 & 1.3 \\
\texttt{t2em           }& 921632   &   4590832 & 4.98  & 1.0 & 1.0 \\
\texttt{tmt\_unsym     }& 917825   &   4584801 & 5.00  & 1.0 & 1.0 \\
\texttt{vas\_stokes\_1M}& 1090664  &  34767207 & 31.88 & 1.0 & 1.3 \\
\texttt{vas\_stokes\_2M}& 2146677  &  65129037 & 30.34 & 1.0 & 1.3 \\
\bottomrule
\end{tabular}
\label{tab:dataset}
\end{table}

We tested several larger-size matrices from the SuiteSparse Matrix Collection~\citep{davis2011university}, the HPCG benchmark~\citep{dongarra2016new}, and the HPGMP benchmark~\citep{yamazaki2022high}, shown in Table~\ref{tab:dataset}.
HPCG is based on the 27-point stencil computation, and the diagonal and off-diagonal elements of the matrices are 26 and -1, respectively.
The matrices from HPGMP are similar to those from HPCG;
the off-diagonal values that represents the connection with forward and backward positions along the z-axis are replaced with $-1 + \beta$ and $-1 - \beta$, respectively ($\beta$ was 0.5 in the experiments).
It should be noted that we applied diagonal scaling to all matrices.
In each test, the right-hand side was a random vector, whose elements were uniformly distributed in the range [0, 1).

We measured the average iteration counts and execution time of three runs.
Starting with an initial guess of a zero vector, each solver ran to reach an approximate solution $\tilde{\vect{x}}$ that satisfies
$$
\frac{\| \vect{b} - \mat{A}\tilde{\vect{x}} \|_2}{\| \vect{b} \|_2} < 1.0 \times 10^{-8}.
$$
When convergence stalled, we terminated \method{} after 300 outermost iterations; that is, \method{} was restarted only three times.
Similarly, FGMRES(64), CG, and BiCGStab were terminated after 19,200 iterations.

\begin{figure*}[tb]
\begin{subfigure}{.95\linewidth}
\includegraphics[width=\linewidth]{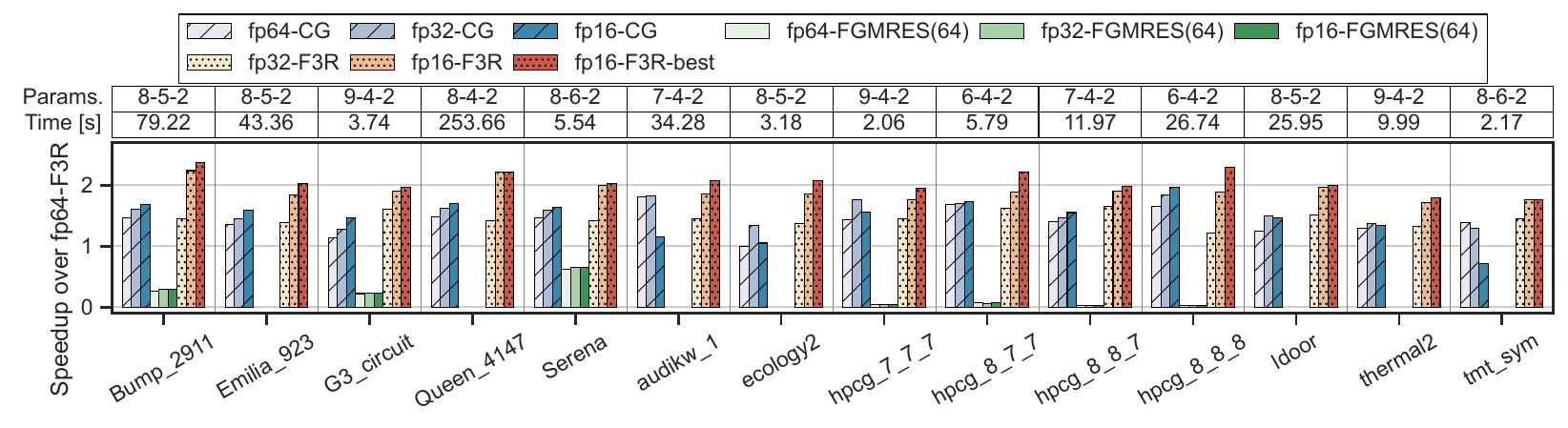}
\caption{Results for symmetric matrices}
\end{subfigure}\\
\begin{subfigure}{.95\linewidth}
\includegraphics[width=\linewidth]{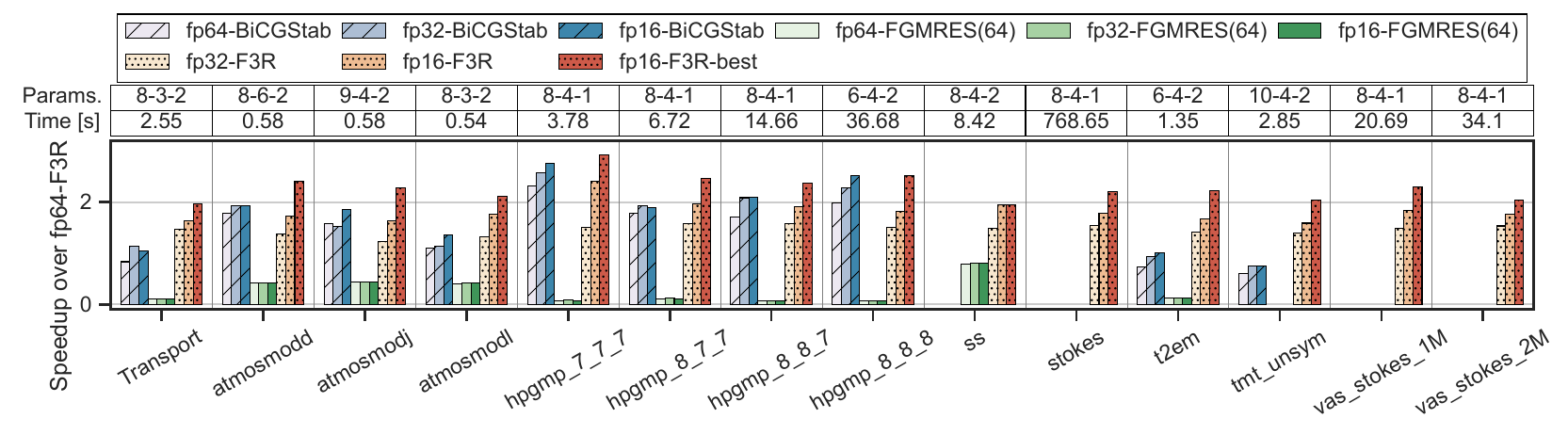}
\caption{Results for non-symmetric matrices}
\end{subfigure}\\
\caption{Performance relative to \methodD{} on the CPU node.
No bar indicates that the convergence failed.
Each table above the plots shows the parameters of \methodHBest{} in $m_2$-$m_3$-$m_4$ at the top row and the execution time of baseline, \methodD{}, in seconds at the bottom row.}
\label{fig:cpu-speedup}
\end{figure*}

\begin{table}[tb]
\small
\caption{Number of invocations of the primary preconditioner $\mat{M}$ until convergence on CPU experiments. Hyphens indicate that the convergence failed.}
\label{tab:cpu-iteration}
\centering
\begin{tabular}{lrrrrr}
\toprule
       & CG or& & fp64- & fp32- & fp16- \\\cmidrule{4-6}
Matrix & BiCGStab & FGMRES(64) & \multicolumn{3}{r}{\method{}} \\
\midrule
\texttt{Bump\_2911} & 4549 & 15833 & 5312 & 5312 & 5248 \\
\texttt{Emilia\_923} & 8696 & - & 11136 & 11072 & 11008 \\
\texttt{G3\_circuit} & 1570 & 2772 & 1664 & 1664 & 1728 \\
\texttt{Queen\_4147} & 6748 & - & 7744 & 7616 & 7616 \\
\texttt{Serena} & 659 & 913 & 832 & 832 & 832 \\
\texttt{audikw\_1} & 2820 & - & 3968 & 4032 & 4032 \\
\texttt{ecology2} & 2421 & - & 2880 & 2816 & 2880 \\
\texttt{hpcg\_7\_7\_7} & 249 & 371 & 384 & 384 & 384 \\
\texttt{hpcg\_8\_7\_7} & 350 & 554 & 512 & 512 & 512 \\
\texttt{hpcg\_8\_8\_7} & 410 & 802 & 512 & 512 & 512 \\
\texttt{hpcg\_8\_8\_8} & 332 & 829 & 512 & 512 & 512 \\
\texttt{ldoor} & 4965 & - & 5824 & 5824 & 5888 \\
\texttt{thermal2} & 3911 & - & 4480 & 4480 & 4544 \\
\texttt{tmt\_sym} & 1761 & - & 2112 & 2112 & 2048 \\
\midrule
\texttt{Transport} & 1060 & 3477 & 896 & 896 & 896 \\
\texttt{atmosmodd} & 214 & 274 & 320 & 320 & 320 \\
\texttt{atmosmodj} & 236 & 254 & 320 & 320 & 320 \\
\texttt{atmosmodl} & 264 & 227 & 256 & 256 & 256 \\
\texttt{hpgmp\_7\_7\_7} & 292 & 416 & 640 & 640 & 512 \\
\texttt{hpgmp\_8\_7\_7} & 368 & 417 & 576 & 576 & 576 \\
\texttt{hpgmp\_8\_8\_7} & 416 & 436 & 640 & 640 & 640 \\
\texttt{hpgmp\_8\_8\_8} & 388 & 512 & 704 & 704 & 768 \\
\texttt{ss} & - & 1246 & 1792 & 1792 & 1664 \\
\texttt{stokes} & - & - & 14976 & 15296 & 15488 \\
\texttt{t2em} & 1770 & 3056 & 1216 & 1216 & 1216 \\
\texttt{tmt\_unsym} & 4374 & - & 2624 & 2624 & 2688 \\
\texttt{vas\_stokes\_1M} & - & - & 4416 & 4416 & 4480 \\
\texttt{vas\_stokes\_2M} & - & - & 3840 & 3840 & 3904 \\
\bottomrule
\end{tabular}
\end{table}

\begin{figure*}[tb]
\begin{subfigure}{.95\linewidth}
\includegraphics[width=\linewidth]{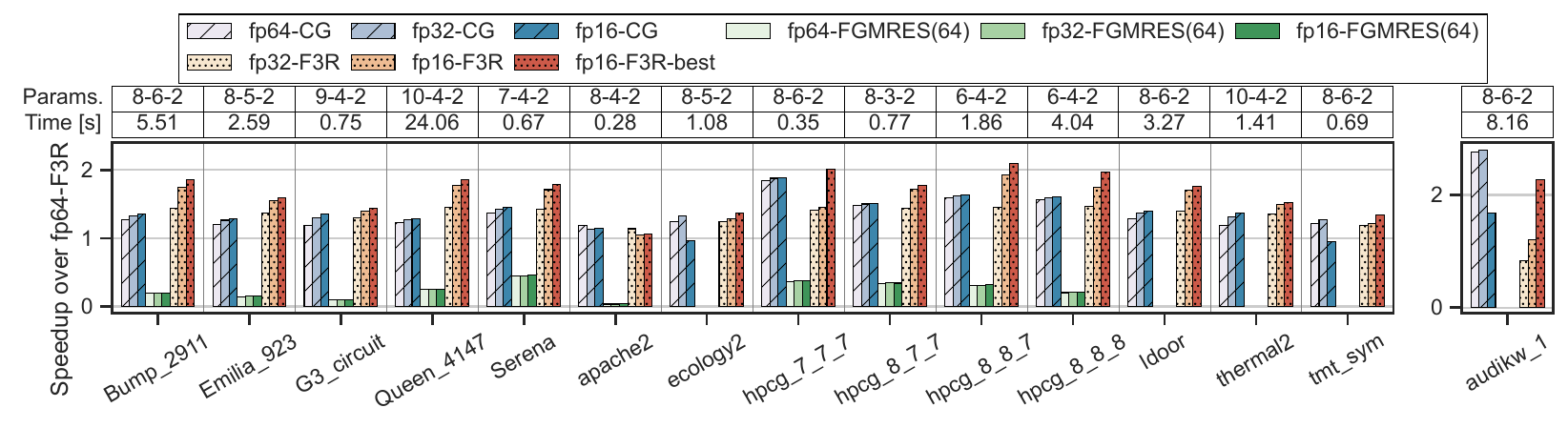}
\caption{Results for symmetric matrices}
\end{subfigure}\\
\begin{subfigure}{.95\linewidth}
\includegraphics[width=\linewidth]{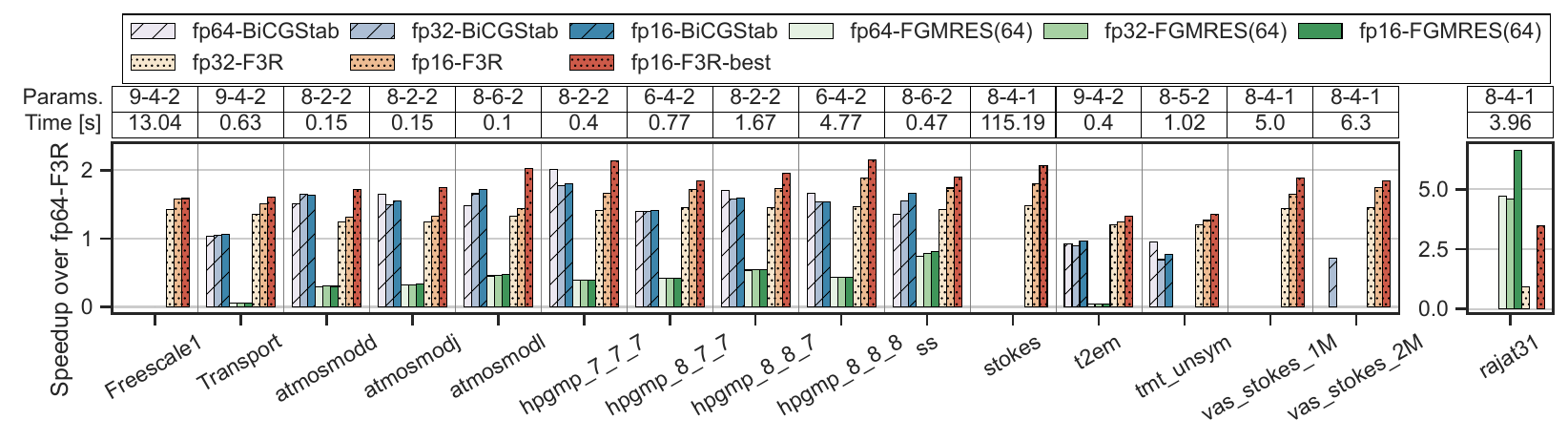}
\caption{Results for non-symmetric matrices}
\end{subfigure}\\
\caption{Performance relative to \methodD{} on the GPU node.
No bar indicates that the convergence failed.
Each table above the plots shows the parameters of \methodHBest{} in $m_2$-$m_3$-$m_4$ at the top row and the execution time of baseline, \methodD{}, in seconds at the bottom row.}
\label{fig:gpu-speedup}
\end{figure*}

\subsection{Experiments on a CPU Node}
First, we conducted experiments on a CPU node of the Camphor~3 supercomputer at Kyoto University.
A block-Jacobi ILU(0) (or IC(0) when symmetric) preconditioner was used as the primary preconditioner $\mat{M}$ for multi-threading.
The number of blocks was 112 (= $56 \times 2$) because the used node has two Intel Sapphire Rapids CPUs with 56 cores.
Additionally, to stabilize block-Jacobi ILU(0), we used a technique that scales diagonal elements of $\mat{A}$ only during the ILU factorization;
specifically, the ILU factorization was applied to a matrix obtained by multiplying the diagonal elements of $\mat{A}$ by a problem-dependent real number, $\alpha_\mathrm{ILU}$.
The value of $\alpha_\mathrm{ILU}$ was determined such that CG or BiCGStab performed better, listed in Table~\ref{tab:dataset}.

The coefficient matrix and preconditioner were stored in the compressed sparse row format.
All vector and matrix-vector operations were multi-threaded in a row-wise manner with 112 threads.
The used compiler was the Intel oneAPI DPC++/C++ Compiler (ver. 2023.2.4) with options \lstinline{-std=c++17} \lstinline{-fast} \lstinline{-qopenmp} \lstinline{-xCORE-AVX512} \lstinline{-mavx512fp16} \lstinline{-qnextgen}.
We also specified \lstinline{numactl} \lstinline{--interleave=all} at runtime to suppress the variability and deterioration of execution performance.

Figure~\ref{fig:cpu-speedup} shows the relative performance, with the upper plot for symmetric matrices and the lower plot for non-symmetric matrices.
The tables accompanying the plots show the average execution time (in seconds) for the baseline \methodD{} (bottom rows) and the values of $m_2$, $m_3$, and $m_4$ for \methodHBest{} in $m_2$-$m_3$-$m_4$ (top rows). 
The results for \lstinline{apache2}, \lstinline{Freescale1}, and \lstinline{rajat31} are omitted because no solver converged.
Table~\ref{tab:cpu-iteration} shows the number of preconditioning steps required for convergence for each of fp64-FGMRES(64), fp64-CG (or fp64-BiCGStab), and the three implementations of \method{}.

First, we analyze the performance balance among the three implementations of \method{}: \methodD{}, \methodF{}, and \methodH{}.
As shown in Table~\ref{tab:cpu-iteration}, there is no significant difference in the convergence rate, regardless of the use of lower-precision arithmetic in \method{}.
Even in the cases where the number of iterations increased, such as \texttt{hpgmp\_8\_8\_8}, the rate was at most 9\%.
Furthermore, in some cases, \methodH{} converged in fewer iterations than \methodD{} and \methodF{}.
These results confirm that using low-precision arithmetic in inner solvers with a limited number of iterations does not slow down convergence.
This compatibility with low-precision arithmetic contributed to the speedups shown in Figure~\ref{fig:cpu-speedup}.
Specifically, \methodF{} outperformed \methodD{} across all tested problems, with an average speedup of about $1.46\times$.
Furthermore, \methodH{} surpassed \methodF{}, achieving speedups over \methodD{} ranging from $1.59 \times$ to $2.42\times$ and over \methodF{} from $1.11\times$ to $1.60\times$.
Unlike previous methods, our approach consistently benefits from \half{} arithmetic, demonstrating the superiority of our approach.

Next, we compare \method{} with other standard Krylov solvers.
In terms of execution performance, \method{} significantly outperformed FGMRES(64) in all tests, regardless of implementation.
Notably, \methodH{} achieved a speedup of up to $58.06\times$ over fp16-FGMRES(64), even when restricted to problems that fp16-FGMRES(64) could solve.
This superiority of \methodH{} stems from nesting FGMRES to reduce the computational cost in the Arnoldi process.
These results suggest that \methodH{} will outperform existing \single{}-\double{} mixed-precision GMRES solvers because their performance is generally expected to be at most twice that of \double{} GMRES.
Regarding the convergence speed in terms of the number of primary preconditioning steps, \method{} also generally outperformed FGMRES(64), converging with up to 74\% fewer preconditioning steps.
Additionally, \method{} was able to solve \lstinline{tmt_unsym} and \lstinline{ecology2} with relatively few iterations, whereas FGMRES(64) failed, reinforcing the idea that nesting is more effective for convergence than restarting.

Compared to CG and BiCGStab, \method{} demonstrated comparable or superior performance.
Specifically, \methodH{} outperformed fp64-CG and fp64-BiCGStab in 14/14 and 12/14 cases, respectively, with speedups ranging from $1.02\times$ to $1.86\times$ for CG and $1.03\times$ to $2.60\times$ for BiCGStab.
Even against fp16-CG and fp16-BiCGStab, \methodH{} performed better in 13/14 and 9/14 cases, respectively.
Table~\ref{tab:cpu-iteration} further shows that the convergence speed of \method{} was comparable to CG and BiCGStab, especially for problems where all solvers required a relatively large number of iterations, such as \texttt{Bump\_2911} and \texttt{Transport}.
On the other hand, for relatively easier problems such as \lstinline{hpcg_8_8_8}, \method{} tended to require noticeably more iterations than CG or BiCGStab.
This may be attributed to the nested structure, which merges the results of inner solvers for segmented Krylov subspaces, rather than searching an overall Krylov subspace.
However, considering that CG is the de facto standard for symmetric positive definite linear systems, the effectiveness of \method{}, which does not rely on symmetry, is particularly noteworthy.
Furthermore, \method{} successfully converged on problems where BiCGStab failed, suggesting that \method{} is more stable due to the monotonically decreasing behavior of FGMRES.

Finally, while \methodH{} performed comparably to \methodHBest{} for some problems, \methodHBest{} achieved the highest performance across all tests.
It was faster than fp16-CG, fp16-BiCGStab, and fp16-FGMRES(64) by up to $2.47\times$, $2.74\times$, and $69.10\times$, respectively.
These results indicate the promise of optimizing the parameters and suggests that \method{} has potential for improvement.
This aspect will be explored in Section~\ref{subsec:restart} with additional experiments.

\subsection{Experiments on a GPU Node}
Next, we present results obtained on a GPU node of the Gardenia supercomputer at Kyoto University.
The node has four NVIDIA A100 GPUs, but the experiments used a single GPU.
We employed the SD-AINV preconditioner~\citep{suzuki2022new}, a simplified version of the standard approximate inverse preconditioner~\citep{benzi1996sparse}, for $\mat{M}$ to exploit the GPU parallelism.
SD-AINV requires only two sparse matrix-vector multiplications (SpMVs) per preconditioning step and is well-suited for GPU implementation.
As with block-Jacobi ILU(0), we scaled the diagonal elements of $\mat{A}$ by a real number $\alpha_\mathrm{AINV}$ during the preconditioner construction.
The value of $\alpha_\mathrm{AINV}$ is listed in Table~\ref{tab:dataset}.

All GPU kernels were implemented in CUDA to execute all computations on the GPU, except for convergence checking.
The NVIDIA HPC Compiler (ver.~23.9), \lstinline{nvc++}, was used, with the main options \lstinline{-std=c++17} \lstinline{-O3} \lstinline{-cuda}.
Unlike the CPU experiments, the SpMV kernels were implemented using the sliced ELLPACK format~\citep{monakov2010automatically} with a chunk size of 32.

The GPU experiments showed results similar to those of the CPU experiments.
As shown in Figure~\ref{fig:gpu-speedup}, \methodF{} outperformed \methodD{}, while \methodH{} was even faster than \methodF{}.
Moreover, all three implementations of \method{} exhibited similar convergence speeds for most problems, indicating that \method{} can effectively incorporate \half{} arithmetic.
While \method{} with the default setting performed well, \methodHBest{} achieved even better results, suggesting room for optimization.
\method{} generally converged faster than FGMRES(64) and performed comparably to CG and BiCGStab.
Specifically, it was faster than FGMRES(64) in all the tests except for \lstinline{rajat31}, achieving speedups ranging from $2.16\times$ to $28.18\times$.
Additionally, \method{} outperformed CG and BiCGStab in many cases; \methodH{} performed better than fp16-CG in 11 out of 15 cases and fp16-BiCGStab in 11 out of 16 cases.
These results confirm that \method{} is also well-suited for GPU implementation and can be effectively combined with the approximate inverse-type preconditioner.

However, two main differences should be noted.
First, in the case of \texttt{rajat31}, \method{} underperformed FGMRES(64), suggesting the existence of a class of problems where nesting FGMRES may degrades the convergence speed; that is, Assumption~\ref{enum:fgmres} does not hold.
Second, the speedup achieved by reducing precision to \single{} and \half{} was moderate, $1.34\times$ and $1.55\times$ on average, compared to $1.46\times$ and $1.87\times$ in the CPU experiments.
This could be attributed to using a different preconditioner and differences in memory bandwidth and cache utilization between CPUs and GPUs.
\method{} requires storing matrix values in \double{}, \single{}, and \half{}, which incurs an overhead and may not have been suitable for the GPU.
To further improve the performance on GPUs, developing a framework that enhances cache efficiency would be necessary.

\begin{figure*}[tb]
\centering
\begin{subfigure}{0.32\linewidth}
\centering
\includegraphics[width=\linewidth]{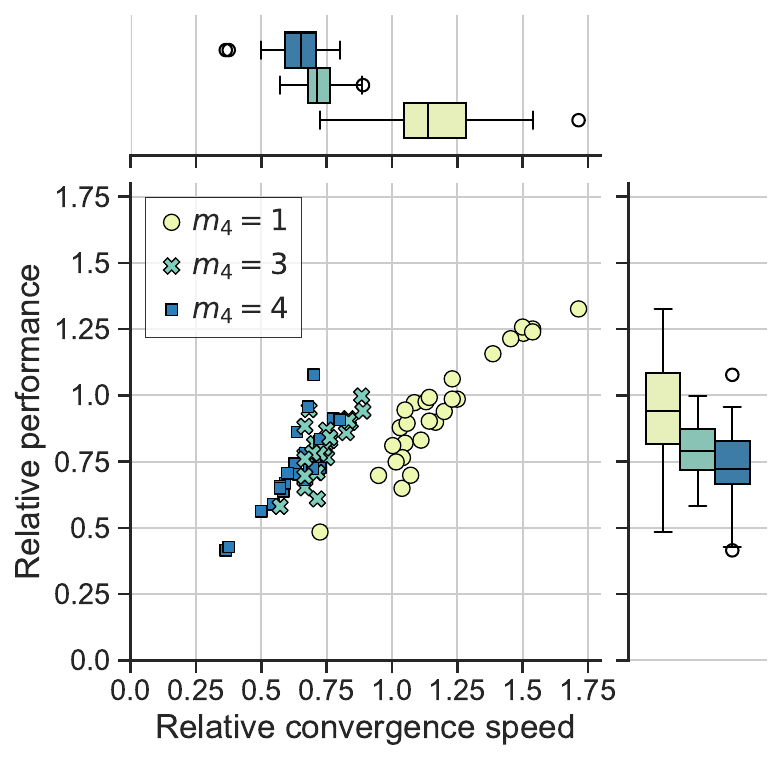}
\caption{Results for $m_4$}
\end{subfigure}\hfill
\begin{subfigure}{0.32\linewidth}
\centering
\includegraphics[width=\linewidth]{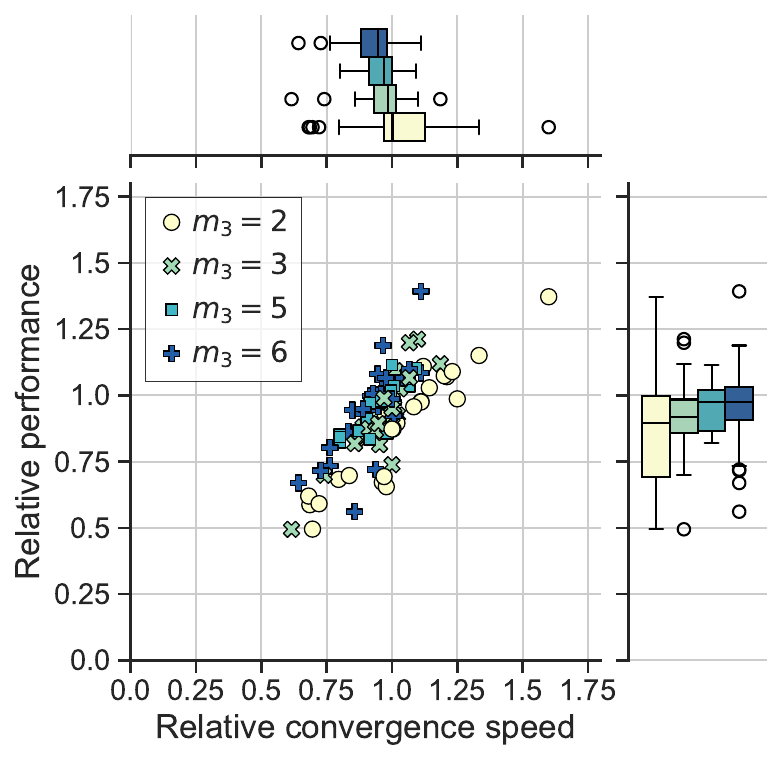}
\caption{Results for $m_3$}
\end{subfigure}\hfill
\begin{subfigure}{0.32\linewidth}
\centering
\includegraphics[width=\linewidth]{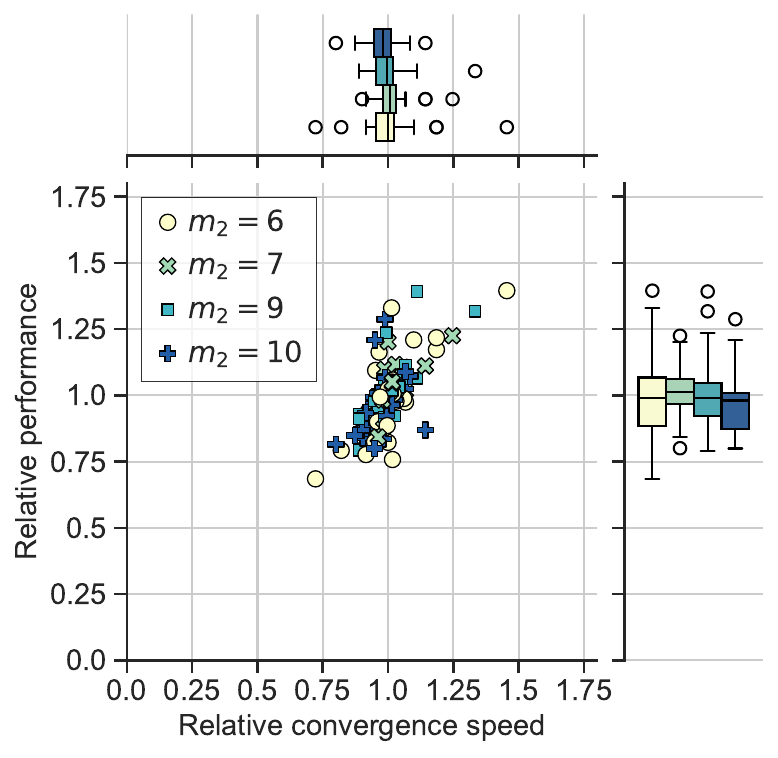}
\caption{Results for $m_2$}
\end{subfigure}
\caption{Results for different values of $m_2$, $m_3$, and $m_4$.
Boxplots on the right and top correspond to the y- and x-axis respectively.
The results are relative to \methodH{} with the default setting ($m_2,m_3,m_4 = 8,4,2$), and the larger, the better in both axes.}
\label{fig:restart}
\end{figure*}

\begin{figure}[htb]
\makeatletter
\if@twocolumn
\includegraphics[width=0.64\linewidth]{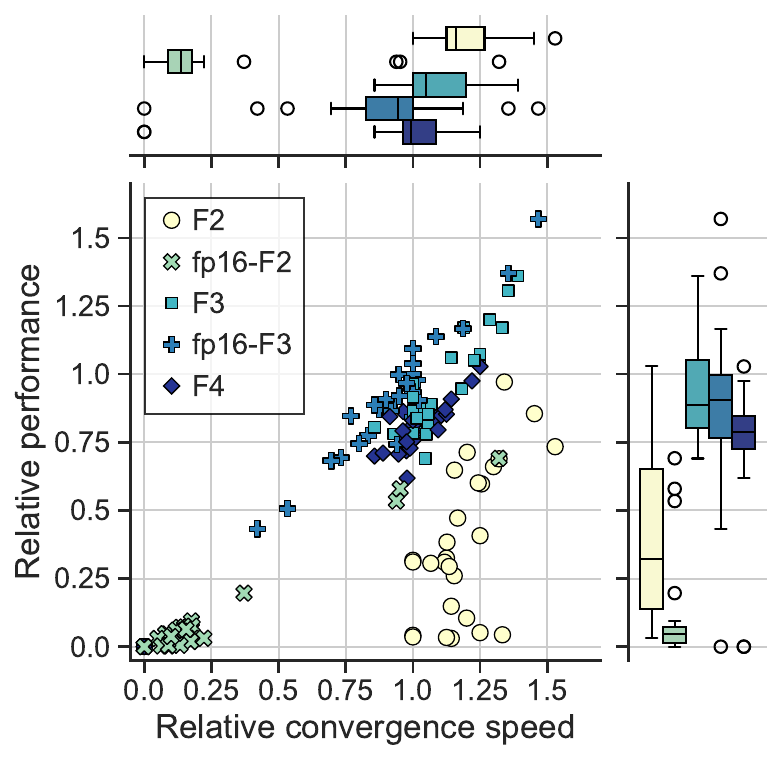}
\else
\centering
\includegraphics[width=0.4\linewidth]{artwork/depthjoint-gray.pdf}
\fi
\makeatother
\caption{Relationship between performance and the depth of nesting.
The results are relative to \methodH{} with the default setting, and the larger, the better.}
\label{fig:depth}
\end{figure}

\begin{figure}[htb]
\makeatletter
\if@twocolumn
\includegraphics[width=0.64\linewidth]{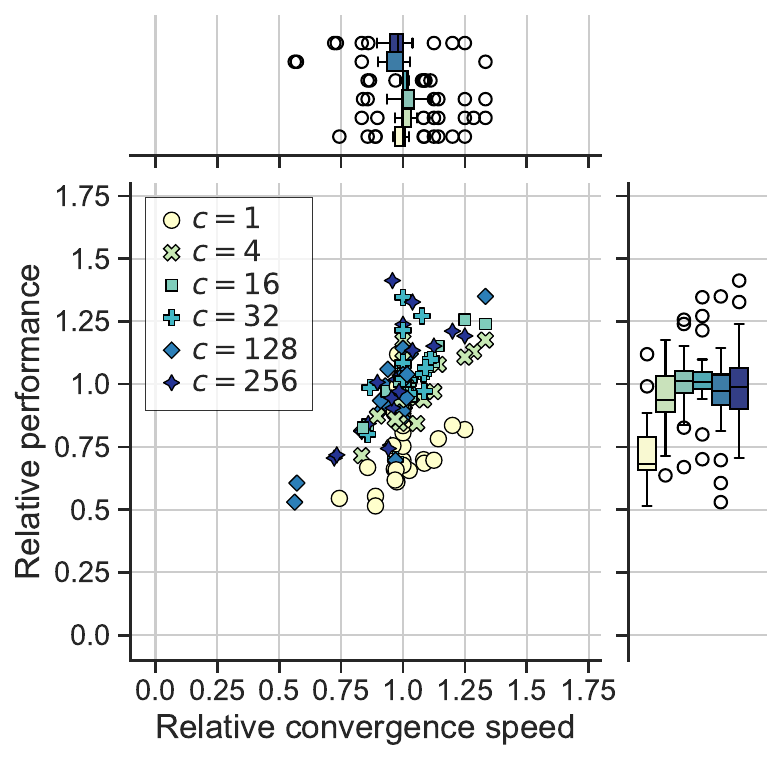}
\else
\centering
\includegraphics[width=0.4\linewidth]{artwork/weight-period-gray.pdf}
\fi
\makeatother
\caption{Performance balance when changing the weight-updating cycle $c$ in the Richardson part.
The results are relative to \methodH{} with $c = 64$.}
\label{fig:weight}
\end{figure}

\begin{figure*}[htb]
\includegraphics[width=.95\linewidth]{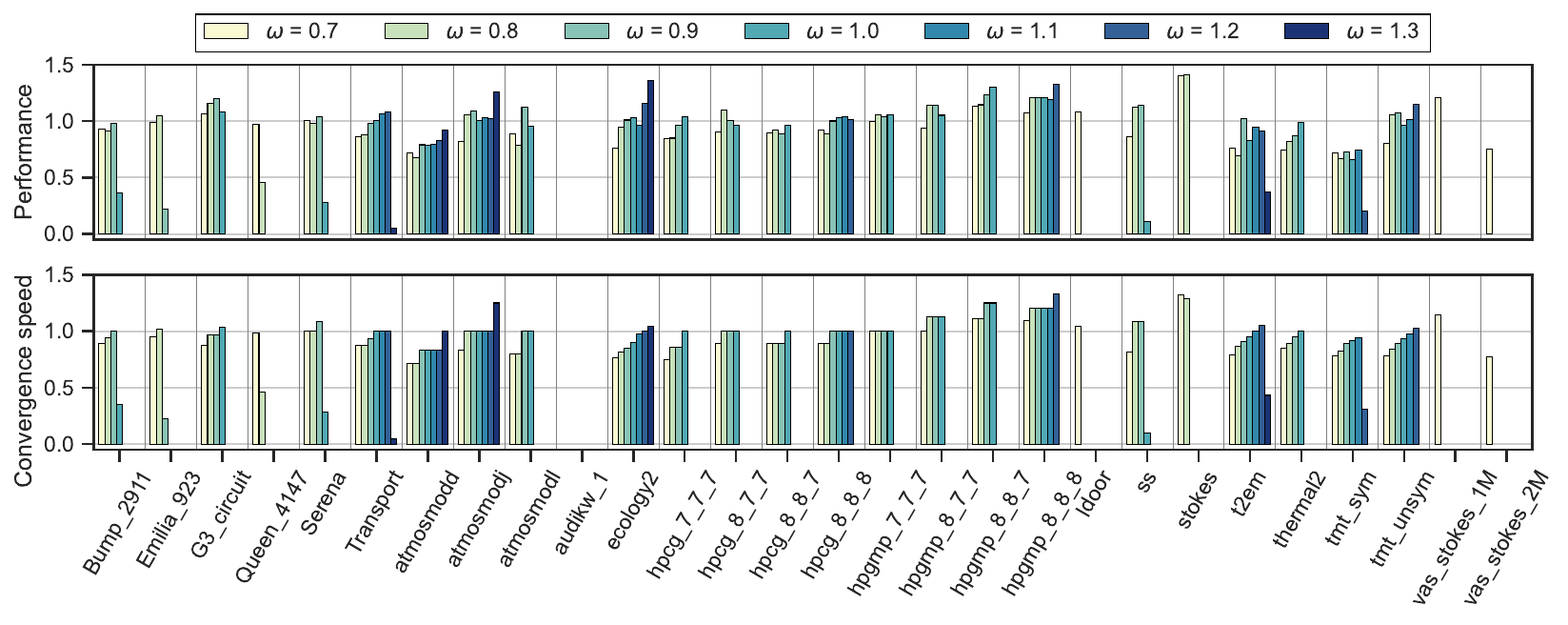}
\caption{\methodH{} with the adaptive strategy vs. \methodH{} with a fixed weight.
No bar indicates that the convergence failed.
The values are relative to \methodH{} with the adaptive strategy.
Values $<1$ indicate that the adaptive strategy was better.}
\label{fig:weight-bar}
\end{figure*}

\section{Experiments on Parameter Setting}
\label{sec:ablation}
In Section~\ref{sec:numerical}, we have focused on the results with the default parameters, $(m_1,m_2,m_3,m_4) = (100,8,4,2)$ and $c = 64$, and demonstrated that these are relatively good but also has room for optimization.
Thus, in this section, we investigate different settings of these parameters with additional experiments on the CPU node.
Furthermore, we discuss the propriety of four-level nesting and using the Richardson iteration in the innermost solver, namely Assumptions~\ref{enum:fgmres} and~\ref{enum:richardson}.

\subsection{On the Number of Inner Iterations}
\label{subsec:restart}
First, we examine the number of iterations of each inner solver, $m_2$, $m_3$, and $m_4$.
Figure~\ref{fig:restart} illustrates trends in the performance of \methodH{} when varying the values of $m_2$, $m_3$, and $m_4$, showing the relationships between execution performance and convergence speed relative to the default setting.
Values $> 1$ on both axes indicate better performance than \methodH{}.
We begin by analyzing the effect of $m_4$, the number of the innermost Richardson iterations.
The results for $m_4 = 3$ and $m_4 = 4$ show that increasing $m_4$ generally reduced the convergence speed rather than improving it.
This degradation directly led to poor execution performance because the innermost Richardson accounts for a certain amount of computation in \methodH{}.
This result suggests that Assumption~\ref{enum:richardson} does not hold for $m = 3, 4$.
In contrast, the setting of $m_4 = 1$ demonstrated better convergence speeds for several problems, though its execution performance was generally lower than that of $m_4 = 2$.
In the case of $m_4 = 1$, the Richardson part, ($\mat{R}_{m_4}$, $\mat{M}$), is equivalent to the standard preconditioner, $\mat{M}$, if the weight is not taken into account, meaning all iterations are performed by FGMRES.
This explains why $m_4 = 2$ often outperforms $m_4 = 1$: a single Richardson iteration involves considerably fewer computations than an FGMRES iteration.

Next, we inspect the effects of $m_3$ and $m_2$.
Compared to $m_4$, both had a smaller effect; that is, no clear relationship was observed between them and the performance of \methodH{}.
Increasing $m_3$ slightly improved the convergence rate; however, since it also increased the computational cost per outermost iteration, the overall execution performance remained relatively unchanged.
Notably, depending on the problem, the performance fluctuated within a range of approximately 50--140\% when comparing $m_3=2$ and $m_3=6$ against the default value $m_3 = 4$.
This suggests that $m_3$, in particular, offers room for optimization, e.g., through adaptive methods.
Similarly, the effect of $m_2$ was small, especially in terms of the variance of the convergence speed; the change was less than 10\% in most tested problems.
The execution performance varied slightly more because of changes in computational cost per outermost iteration, but its effect was still minor compared to $m_4$ and $m_3$.

Consequently, while the default choice $m_4 = 2$ is reasonable, finding optimal values of $m_2$ and $m_3$ is more challenging as they depend strongly on the problem.
In particular, $m_3$ has a comparatively large effect, which implies the importance of techniques to adaptively control it for further improved performance of \method{}.

\newcommand{\depthII}{F2}
\newcommand{\depthIII}{F3}
\newcommand{\depthIV}{F4}
\newcommand{\depthIIH}{fp16-F2}
\newcommand{\depthIIIH}{fp16-F3}
\newcommand{\depthIVH}{fp16-F4}

\begin{table}
\caption{Additional solvers used in Section~\ref{subsec:depth}}
\label{tab:depth}
\small
\centering
\begin{tabular}{lllrrr}
\toprule
     &                &  & \multicolumn{3}{c}{Precision} \\\cmidrule{4-6}
Solver & Tuple notation & Parts & $\mat{A}$ & Vectors & $\mat{M}$ \\
\midrule
\depthII{}   & ($F_{100}$, $F_{64}$, $M$) & $F_{100}$ & \double{} & \double{} & - \\
             &                            & $F_{64}$  & \single{} & \single{} & \half{} \\ \midrule
\depthIIH{}  & ($F_{100}$, $F_{64}$, $M$) & $F_{100}$ & \double{} & \double{} & - \\
             &                            & $F_{64}$  & \half{}   & \half{}   & \half{} \\ \midrule
\depthIII{}  & ($F_{100}$, $F_8$, $F_8$, $M$) & $F_{100}$ & \double{} & \double{} & - \\
             &                                & $F_{8}$ (outer)  & \single{} & \single{} & - \\
             &                                & $F_{8}$ (inner)   & \half{}   & \single{} & \half{} \\ \midrule
\depthIIIH{} & ($F_{100}$, $F_8$, $F_8$, $M$) & $F_{100}$ & \double{} & \double{} & - \\
             &                                & $F_{8}$ (outer)  & \single{} & \single{} & - \\
             &                                & $F_{8}$ (inner)   & \half{}   & \half{} & \half{} \\ \midrule
\depthIV{}   & ($F_{100}$, $F_8$, $F_4$, $F_2$, $M$) & $F_{100}$ & \double{} & \double{} & - \\
             &                                       & $F_{8}$   & \single{} & \single{} & - \\
             &                                       & $F_{4}$   & \half{} & \single{}   & - \\
             &                                       & $F_{2}$   & \half{} & \half{}     & \half{} \\
\bottomrule
\end{tabular}
\end{table}

\subsection{On the Nesting Depth} \label{subsec:depth}
Next, we investigate the depth of nesting by showing results for five additional solvers in Figure~\ref{fig:depth}.
The characteristics of these solvers are summarized in Table~\ref{tab:depth}.
Briefly, \depthII{} and \depthIII{} are two- and three-level nested FGMRES, respectively.
Each maintains the same precision as \methodH{} at corresponding depth levels.
\depthIIH{} and \depthIIIH{} are modified versions of \depthII{} and \depthIII{}, respectively, where the innermost solver is implemented in \half{}.
\depthIV{} is largely the same as \methodH{} but replaces the innermost Richardson part with FMGRES.

Before considering the depth of nesting, we compare \methodH{} and \depthIV{} to verify Assumption~\ref{enum:richardson}.
The figure shows that the convergence rates of \depthIV{} and \methodH{} were similar, with a relatively small quartile deviation (all plots are close to 1.0), thereby validating Assumption~\ref{enum:richardson} with $m_4 = 2$.
Furthermore, despite the similar convergence speed, \methodH{} outperformed \depthIV{} for all problems except for \lstinline{t2em}.
This demonstrates that using Richardson is advantageous in terms of computational cost because it skips the Arnoldi process.

Next, we examine \depthII{} and \depthIIH{} which can be regarded as \double{}-\single{} and \double{}-\half{} inner-outer FGMRES, respectively.
\depthII{} was generally better than \methodH{} in the convergence speed, but the improvement was limited, falling below 1.25 times in about 75\% of the test cases.
This suggests that Assumption~\ref{enum:fgmres} is not entirely implausible;
we can replace FGMRES with a nested one without significantly degrading the convergence.
In execution performance, \depthII{} underperformed \methodH{} because of its higher computational cost in the Arnoldi process, as mentioned in Section~\ref{subsec:guide-to-nesting}, indicating that nesting FGMRES has certain merit.
Moreover, \depthIIH{} significantly reduced the convergence speed and, in turn, execution performance.
This reaffirms that using \half{} for all parts of inner (flexible) GMRES with long iterations causes a precision overflow.

Finally, we evaluate the three-level nested solvers, \depthIII{} and \depthIIIH{}.
To the best of our knowledge, three-level nested FGMRES with \single{} and \half{} has not been previously reported.
For more than half of the tested problems, the proposed solver, \methodH{}, outperformed both \depthIII{} and \depthIIIH{} in execution performance, underscoring the effectiveness of \methodH{}.
Moreover, the results also offer two key observations.
First, \depthIII{} and \depthIV{} had similar convergence speeds for many problems.
This indicates that approximating ($\mat{F}_8$, $\mat{M}$) with ($\mat{F}_4$, $\mat{F}_2$, $\mat{M}$) did not significantly degrade convergence, supporting the validity of Assumption~\ref{enum:fgmres}.
However, for several problems, such as $\texttt{stokes}$ and $\texttt{vas\_stokes\_2M}$, \depthIII{} converged faster, resulting in better execution performance even compared to \methodH{}.
These problems are considered relatively difficult, as neither BiCGStab nor FGMRES(64) could solve them.
This implies that, for difficult problems, searching a large subspace is more effective than combining results from small subspaces via the nested Krylov approach.
In these cases, Assumption~\ref{enum:fgmres} does not hold, and we must carefully consider the number of inner iterations for better performance.
Second, \depthIIIH{} showed slower convergence than \depthIII{}, indicating that eight iterations are still excessive to perform FGMRES in \half{} without exceeding precision limits.
Therefore, limiting the number of \half{} iterations appears to be a practical strategy, and four-level nesting in \methodH{} is a reasonable choice to balance convergence, execution performance, and the effective use of \half{}.

\subsection{On the Weight in Richardson}
Finally, we discuss the weight in the Richardson part.
Figure~\ref{fig:weight} compares six different values of the updating cycle, $c$, in the adaptive weight-updating technique.
No clear trend emerges compared to the other factors discussed above, except for the case of $c = 1$.
In this case, the convergence speed was not improved significantly, despite many additional computations imposed on \method{}.
As a result, the default setting $c = 64$ performed similar to the case of $c = 1$ for most cases.
For other settings, performance varied depending on the problem.
Looking more closely, a slight decrease in convergence is observed as $c$ increases; however, execution performance does not necessarily degrade because of the fewer additional computations due to the large $c$.

Figure~\ref{fig:weight-bar} compares our adaptive approach with a static one, where the weight is manually set.
In some cases, such as $\omega=1.3$ for \lstinline{ecology2}, the static approach performed better.
However, the adaptive approach was one of the best in most cases.
Importantly, the performance of the static approach is highly sensitive to the setting.
For instance, the static approach failed on \lstinline{audikw_1} for all weight values, whereas the dynamic approach succeeded.
This stability against the parameter makes our adaptive technique a viable and effective option.

\subsection{Summary}
In this section, we have analyzed \method{} from three perspectives: the number of inner iterations, the nesting depth, and the weight in Richardson.
Experimental results confirmed the validity of both Assumptions~\ref{enum:richardson} and~\ref{enum:fgmres} within the context of our approach.
However, Assumption~\ref{enum:fgmres} occasionally did not hold, particularly for relatively difficult problems.
Additionally, the results showed that the default parameters $(m_1,m_2,m_3,m_4) = (100, 8,4,2)$ effectively balanced performance, convergence, and the use of \half{}.
Notably, the value of $m_3$ had a large effect on the performance of \method{}, though its optimal value varied depending on problem.
Since identifying it is challenging, adaptive approaches, such as dynamically controlling it based on the magnitude of the residuals, would be beneficial.
Finally, it was confirmed that the adaptive weight-updating technique in Algorithm~\ref{alg:richardson} provided greater stability to the parameter than using a single, fixed weight.

\section{Conclusions} \label{sec:conclution}
We proposed a solver, \method{}, based on a novel mixed-precision approach that integrates FGMRES and Richardson solvers in the nested Krylov structure.
Numerical experiments confirmed that this approach is suited for using \half{} arithmetic without increasing the total number of iterations.
Additionally, the results showed that \method{} is more stable than BiCGStab and achieves comparable or superior performance to CG, BiCGstab, and restarted FGMRES.
This is because \method{} requires relatively fewer computations in the Arnoldi process compared to restarted FGMRES.
Experiments presented in Section~\ref{sec:ablation} indicated that the parameters used in Section~\ref{sec:numerical} are reasonable, though there is still potential for improvement, especially in the number of iterations performed by the middle solvers.

Future direction includes the following three.
The first is conducting more comprehensive evaluations.
Since \method{} relies solely on FGMRES and Richardson, it can apply asynchronous preconditioners, which could also be advantageous.
Moreover, other Krylov methods also accept various approaches to reducing precision beyond simply using lower-precision preconditioners, although these would need parameter optimization.
These methods may change the performance balance.
The second direction is introducing dynamic techniques for terminating inner iterations.
For instance, a strategy based on magnitude of inner residuals could be beneficial.
The third direction is extending \method{} for distributed systems.
Notably, the innermost Richardson part does not incorporate a reduction process, which typically involves collective communication.
Furthermore, the inner FGMRES solvers perform only a limited number of iterations per parent iteration, making them well-suited for communication-avoiding implementations~\citep{hoemmen2010communication, yamazaki2024two}.

\section*{Acknowledgement}
This work was supported by JSPS KAKENHI Grant Number JP24KJ0266.

\bibliographystyle{ACM-Reference-Format}
\bibliography{fxr}

\end{document}